\documentclass{article}

\usepackage{arxiv}

\usepackage[utf8]{inputenc} 
\usepackage[T1]{fontenc}    
\usepackage{hyperref}       
\usepackage{url}            
\usepackage{booktabs}       
\usepackage{amsfonts}       
\usepackage{nicefrac}       
\usepackage{microtype}      
\usepackage{graphicx}
\usepackage{doi}

\usepackage{graphicx}
\usepackage{caption}
\usepackage{subcaption}
\usepackage{bm}
\usepackage{physics}
\usepackage{amsmath}
\usepackage{mathtools}
\usepackage{xcolor}

\newcommand{\rr}{\mathbf{r}}
\newcommand{\pp}{\mathbf{p}}

\title{Globally time-reversible fluid simulations with smoothed particle hydrodynamics}

\author{Ond{\v r}ej Kincl\\
Mathematical Institute\\
Faculty of Mathematics and Physics\\
Charles University\\
Sokolovská 49/83, Prague, 186 75\\
Czech Republic\\
	\texttt{ondrej.kincl.6@gmail.com} \\
	\And
	\href{https://orcid.org/0000-0003-0605-6737}{\includegraphics[scale=0.06]{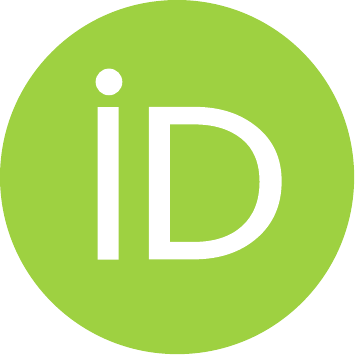}\hspace{1mm}Michal Pavelka} \\
Mathematical Institute\\
Faculty of Mathematics and Physics\\
Charles University\\
Sokolovská 49/83, Prague, 186 75\\
Czech Republic
}

\renewcommand{\shorttitle}{Globally Reversible SPH}

\hypersetup{
pdftitle={Globally time-reversible fluid simulations with smoothed particle hydrodynamics},
pdfauthor={Ond{\v r}ej Kincl, Michal Pavelka},
pdfkeywords={smoothed particle hydrodynamics, symplectic, reversibility}
}

\begin{document}
\maketitle

\begin{abstract}
			This paper describes an energy-preserving and globally time-reversible code for weakly compressible smoothed particle hydrodynamics (SPH). We do not add any additional dynamics to the Monaghan's original SPH scheme at the level of ordinary differential equation, but we show how to discretize the equations by using a corrected expression for density and by invoking a symplectic integrator. Moreover, to achieve the global-in-time reversibility, we have to correct the initial state, implement a conservative fluid-wall interaction, and use the fixed-point arithmetic. Although the numerical scheme is reversible globally in time (solvable backwards in time while recovering the initial conditions), we observe thermalization of the particle velocities and growth of the Boltzmann entropy. In other words, when we do not see all the possible details, as in the Boltzmann entropy, which depends only on the one-particle distribution function, we observe the emergence of the second law of thermodynamics (irreversible behavior) from purely reversible dynamics. 
\end{abstract}

\keywords{SPH \and reversibility \and symplectic integrator \and stability \and entropy}

\section{Introduction}
The original SPH scheme, now called weakly compressible SPH (WCSPH), was developed in 1960s by R.Gingold and J. Monaghan \cite{sph}. The scheme used the leap-frog time integration and relied on an artificial term stabilizing the fluid in the presence of shocks. Many other stabilization strategies were conceived, including $\delta$-SPH \cite{antuono2010} and methods based on Riemann solvers \cite{Zhang2017}. These approaches often provide results in reasonable agreement with experimental data or other numerical methods. Unfortunately, they often come at the cost of increased complexity and necessity to fine-tune additional numerical parameters. Failure in such fine-tuning can result in unrealistic dissipation, or large energy oscillations unless a convenient limiter is used \cite{meng2021}.
	
	There are essentially two different ways, how pressure can be calculated in WCSPH \cite{violeau}. A commonly used approach is to update the density iteratively, using the velocities obtained from the discretized balance of momentum. Alternatively, one can compute density directly as an SPH interpolation of particle masses. The latter approach is less common, but it leads to a symplectic structure of the system of SPH differential equations, which makes the WCSPH discretization completely stable even in the absence of viscosity or stabilization. Moreover, the WCSPH numerical scheme becomes globally reversible in time when also the fixed-point arithmetic is used. For instance, the breaking dam simulation can be reversed to recover the initial conditions. The purpose of this paper is to show key advantages of the symplectic WCSPH compared to the standard approach.
	
	However, as pointed out by D. Violeau, symplectic WCSPH shows unphysical behavior at free surfaces, where constant functions are not correctly interpolated. We suggest two ways how to correct this issue. The first method modifies the density by an appropriate integration constant, which is easy to implement. The second and more complex methods renormalizes the initial state by solving a certain nonlinear problem using Newton's method. We also demonstrate that a careful implementation of the symplectic WCSPH scheme, using the corrected treatment of free surfaces and the fixed-point (as opposed to the more usual floating-point) arithmetic leads to bit-perfect global time reversibility. This can be considered an additional symmetry-preserving feature alongside the conservation of energy, mass, momentum and angular momentum.  
	
	Although the resulting code is energy-preserving and reversible, it can still be considered dissipative in the sense that the Boltzmann entropy of the SPH particles grows. In other words, we observe the emergence of the second law of thermodynamics from purely reversible dynamics.
	
	We demonstrate the various versions of WCSPH on dam-break and Gresho vortex benchmarks.
	
	All SPH codes used to make this paper are freely available within the new SmoothedParticles.jl package \cite{Kincl_SmoothedParticles_jl} written in the Julia programming language \cite{julia2017}.
	
	\section{WCSPH for inviscid fluids}
	The standard weakly compressible spatial semi-discretization of Euler fluid equations in a uniform gravitational field is \cite{monaghan2005}:
	\begin{subequations}\label{eq.WCSPH}
		\begin{align}
			\partial_{t} \rho_a &= \sum_b m_b (\bm{u}_a - \bm{u}_b) w'(r_{ab}) \bm{e}_{ab}, \label{boma} \\
			\partial_{t} \bm{u}_a  &= \sum_b m_b \left(\frac{p_a}{\rho_a^2} + \frac{p_b}{\rho_b^2}\right) w'(r_{ab}) \bm{e}_{ab} - g \bm{\hat{z}}, \label{bomo} \\
			\partial_{t}\bm{r}_a &= \bm{u}_a, \label{move}
		\end{align}
	\end{subequations}
	where the pressure is a function of density via the baratropic formula,
	\begin{equation}
		p(\rho) = \frac{c^2\rho_0}{7} \left[ \left(\frac{\rho}{\rho_0}\right)^7 - 1\right].
	\end{equation}
	Constant $\rho_0$ stands for the referential density, $g$ is the gravitational acceleration, and $c$ is the numerical speed of sound (typically chosen as ten times the characteristic flow speed $U_\text{char}$). It is well known that the system of ordinary differential equations \eqref{eq.WCSPH} conserves the total energy in the form
	\begin{equation}
		\mathcal{H} = \sum_a m_a\left( \frac{v_a^2}{2} + 
		\epsilon_a +  g z\right),
		\label{boe}
	\end{equation}
	where
	\begin{equation*}
		\epsilon_a = \frac{m_a c^2}{42} \left[ \left( \frac{\rho}{\rho_0}\right)^6 + \frac{6\rho_0}{\rho} \right].
	\end{equation*}
	Moreover, for $g = 0$, the momentum
	$$ \bm{\mathcal{M}} = \sum_{a} m_a \bm{u}_a, $$
	and angular momentum
	$$ \bm{\mathcal{L}} = \sum_{a} m_a \bm{r}_a \times \bm{u}_a$$
	are conserved as well. Naturally, gravitational force can deliver some momentum to the fluid, while the equal opposite reaction of the fluid on Earth is considered grounded.
	
	
	The WCSPH equations \eqref{eq.WCSPH} represent a Hamiltonian system, which can be checked for instance using program \cite{kroeger2010}. But it is not a symplectic system because the mass density is treated as a state variable. If, on the other hand, the density would always be calculated from the current positions of the particles, then only the positions and momenta would be state variables and the system of equations would be symplectic. Symplecticity would make it easier to choose a geometric structure-preserving numerical integrator, for instance using the Verlet scheme.
	
	A symplectic form of WCSPH is obtained by solving the ordinary differential equation for density \eqref{boma},
	\begin{equation}
		\rho_a = \sum_{b} m_b w(r_{ab}),
		\label{boma-closed-0}
	\end{equation}
	instead of updating $\rho_a$ iteratively by \ref{boma}. However, this closed expression for density \eqref{boma-closed-0} is rarely used in practice because it leads to questionable behavior near free surfaces \cite{violeau}. This happens because free boundary particles are under-occupied and, according to ($\ref{boma-closed-0}$), they have $\rho \ll \rho_0$. Therefore, they are equipped with large internal energy and, affected by its negative gradient, they start to vibrate. A remedy is provided in the following section, which leads to the possibility to use the symplectic form of WCSPH.
	
	\section{Treatment of free surfaces}
	\subsection{Solving the equation for $\rho_a$ with an appropriate integration constant}
	The issue of unphysical behavior at free surfaces can be corrected by choosing an appropriate integration constant when solving the differential equation for $\rho_a$,
	\begin{equation}
		\rho_a = \sum_{b} m_b w(r_{ab}) + C_a,
		\label{boma-closed}
	\end{equation}
	where
	\begin{equation*}
		C_a = \rho_0 - \left.\rho_a \right|_{t=0}. 
	\end{equation*}
	Indeed, since WCSPH based on \eqref{boma} typically assumes zero internal forces at $t=0$, it is actually \eqref{boma-closed} and not \eqref{boma-closed-0} that is equivalent to the standard WCSPH formulation, since \eqref{boma-closed-0} leads to non-zero density gradients near the boundary. Although schemes based on \eqref{boma} and \eqref{boma-closed} are completely equivalent at the ODE level, they become very much different once time discretization is considered. It appears that \eqref{boma} is currently used by most researchers, despite the advantages of Equation \eqref{boma-closed}:
	\begin{enumerate}
		\item Inferring $\rho$ from positions using Equation \eqref{boma-closed} does not accumulate the density error.
		\item The system of SPH evolution equations becomes symplectic (and not only Hamiltonian).
		\item Therefore, symplectic structure-preserving integrators can be employed, which leads to simulations globally reversible in time.
	\end{enumerate}
	
	However, even if we use the correct closed formula for density \eqref{boma-closed} and use a symplectic integrator, are facing another obstacle.
	For problems involving free surfaces, integration constant $C_a$ is typically much bigger for boundary particles and, therefore, whenever such particle submerges into interior of the fluid body, it generates a small void pocket around itself. This unwanted numerical artifact is removed in the following section.
	
	\subsection{Initial state correction (ISC)} \label{sec.ISC}
	The numerical artifact that particles that are originally at the boundary are treated differently than the rest of the particles, even if the submerge into the in the interior of the fluid, is rooted in the difference in the $C_a$ constants at the beginning. Therefore, w´e shall remove that artifact by correcting the initial state. 
	
	A solution is to correct the initial positions by a vector field $\delta \bm{x}$ such that
	$$ \rho_0 = \sum_{b} m_b w\left(\big|\bm{x}_a + {\delta \bm{x}_a} - \bm{x}_b\big|\right).$$
	These non-linear equations can be solved by a strategy based on the Newton's method, details of which follow. Linearization with respect to $\delta \bm{x}_a$ leads to
	$$ \rho_0 - \rho_a \doteq \sum_b m_b (\delta \bm{x}_a - \delta \bm{x}_b) \cdot w'(r_{ab}) \bm{e}_{ab}.$$
	On the left hand side, we identify the discrete divergence operator $\tilde{D}^1$ 
	\begin{subequations}
		\label{init renorm}
		\begin{equation}
			-\rho_a \tilde{D}^1_a \delta \bm{x} \doteq  \rho_0 - \rho_a,
			\label{init renorm eq 1}
		\end{equation}
		see \cite{violeau} for details. 
		This linear problem for unknowns $\delta \bm{x}$ is underdetermined, so shall be restricted to $\delta \bm{x}$ in the form of a discrete gradient
		\begin{equation*}
			\delta \bm{x}_a = \sum_b m_b \left(\frac{\phi_a}{\rho_a^2} + \frac{\phi_b}{\rho_b^2}\right) w'(r_{ab}) \bm{e}_{ab},
		\end{equation*} 
		which can be also rewritten as
		\begin{equation}
			\delta \bm{x}_a + \frac{1}{\rho_a} G^1_a \phi = 0.
			\label{init renorm eq 2}
		\end{equation}
	\end{subequations}
	Sparse linear system \eqref{init renorm} for $\delta \bm{x}$ and $\phi$ has to be solved iteratively, updating $\bm{x}$ and reevaluating $\rho$ according to \eqref{boma-closed-0}, until $\rho \doteq \rho_0$ with satisfactory precision. 
	This initial state correction (ISC) procedure has been tested in our numerical examples, and the procedure converges quadratically except for problems with too much symmetry (for instance rectangle with particles arranged in square grid nearly always diverges)\footnote{For these cases, we recommend to initially disrupt particle positions by a random noise.}. The whole procedure may be then summarized as follows:
	\begin{figure}[ht!]
		\centering
		\begin{subfigure}{0.45\linewidth}
			\includegraphics[width=\linewidth]{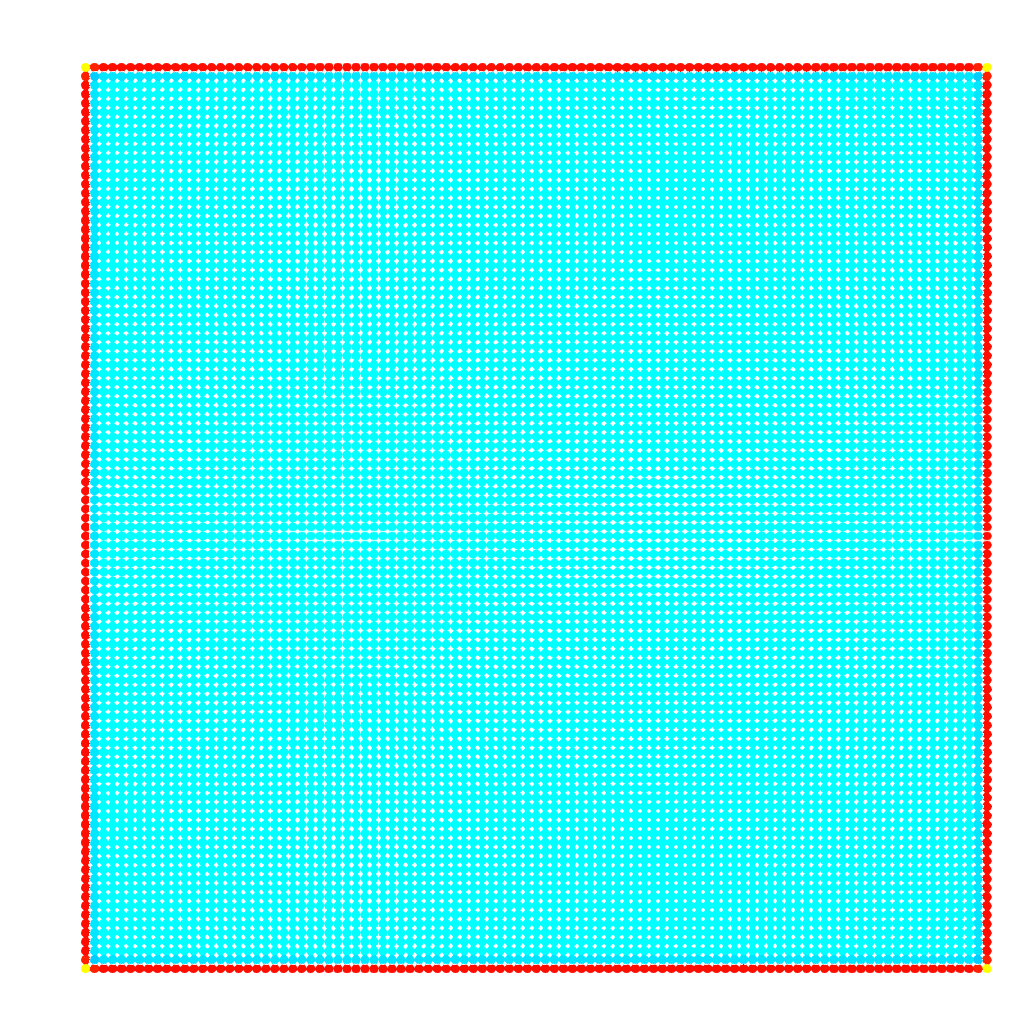}
			\subcaption{before ISC}
		\end{subfigure}
		\begin{subfigure}{0.45\linewidth}
			\includegraphics[width=\linewidth]{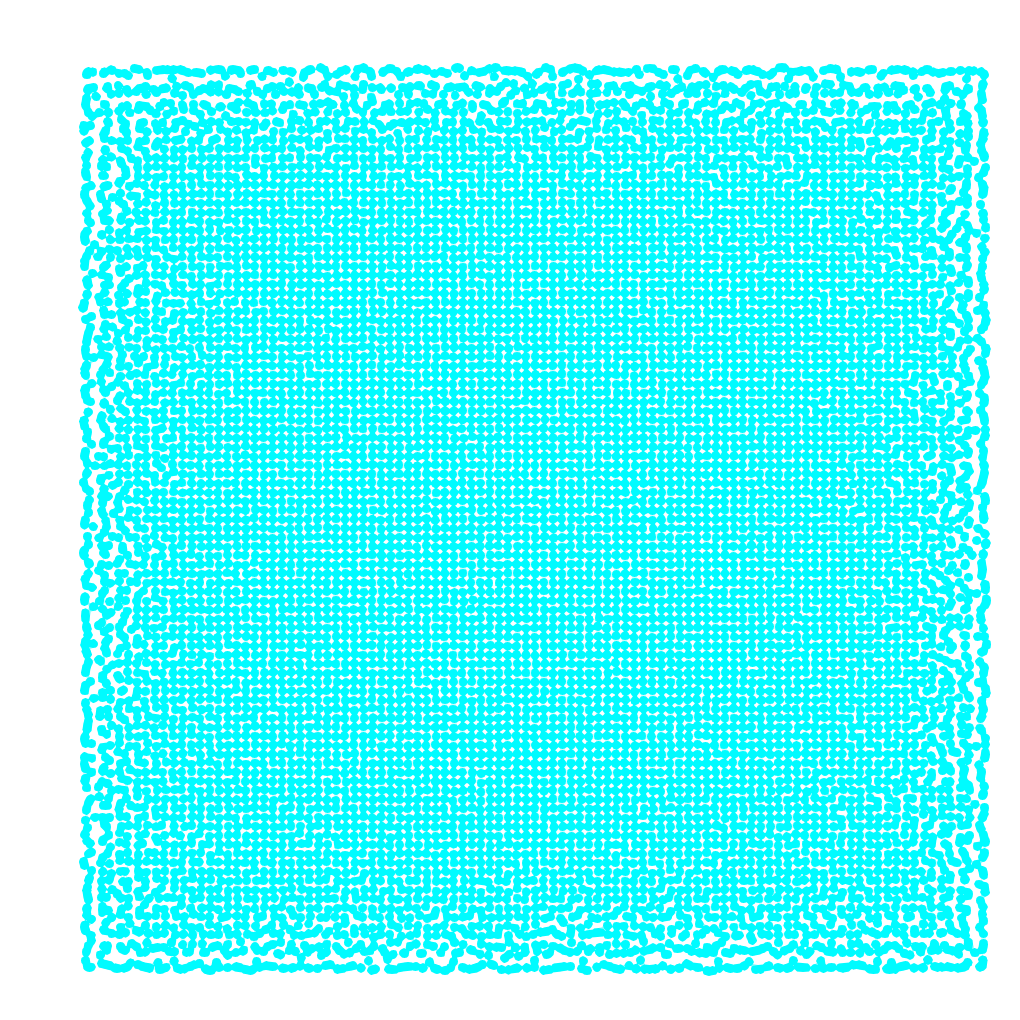}
			\subcaption{after ISC}
		\end{subfigure}
		\caption{Applying ISC to particles in square arranged regularly. Density \eqref{boma-closed-0} is shown in color (red to blue). Initially, it varies due to boundary effect but becomes everywhere constant after 7 iterations of ISC (up to relative error less than $10^{-10}$). }
	\end{figure}
	\begin{center}
		\fcolorbox{black}{white}{\parbox{0.75\textwidth}{
				\centering{\textbf{Initial state correction (ISC)}}
				\begin{enumerate}
					\item set $|\delta \bm{x}_a| < \frac{\delta r}{10}$ randomly for all $a$
					\item repeat until sufficient precision is acquired:
					\item \qquad $\bm{x}_a = \bm{x}_a + \delta \bm{x}_a$ for all $a$
					\item \qquad recompute neighbor list
					\item \qquad find $\rho_a = \sum_{b} m_b w(r_{ab})$ for all $a$
					\item \qquad solve $N^{d+1} \times N^{d+1}$ linear problem
					$$  \begin{bmatrix}
						I & \frac{1}{\rho} G^1\\
						-\rho \tilde{D}^1 & 0
					\end{bmatrix} 
					\begin{bmatrix}
						\delta \bm{x}  \\
						\phi 
					\end{bmatrix}
					= \begin{bmatrix}
						0 \\
						\rho_0 - \rho
					\end{bmatrix}$$
		\end{enumerate}}}
	\end{center}
	
	This computation is expensive, but is performed only once per simulation. In place of the linear system with $2\times 2$ block matrix, one could also solve directly for $\phi$:
	$$ -\rho \tilde{D}^1 \, \frac{1}{\rho} G^1 \phi = \rho_0 - \rho.$$
	However, the product $A = -\rho \tilde{D}^1 \, \frac{1}{\rho} G^1$ leads to a matrix which is much more dense. A matrix-free iterative method could be used without evaluating the product directly, is a future possibility.
	
	Let us note that after ISC, all constant fields will be correctly interpolated, at least initially. Weak compressibility helps to preserve this trait approximately in a way that does not accumulate error provided \eqref{boma-closed-0} is used. In this sense, ISC can be considered an alternative to (zeroth order) operator renormalization with the benefit that it does not break the symmetry of smoothing kernel. Unfortunately, this technique has two drawbacks. Firstly, it may slightly deform geometry in an unwanted way. Secondly, adapting it to a concrete setting may be difficult --- for instance, when an anti-clump term is included. 
	
	Both symplecticity of the system of equations and ISC improve the numerical solutions, but in order to obtain a globally reversible behavior we have to tackle a further problem caused at the walls of the simulation box, as in the following section.
	
	\subsection{Wall modelling}
	In WCSPH, walls can be modeled by a layer of dummy particles which behave as if they belonged to the fluid, except that their positions and velocities are fixed. In this approach, it is hoped that the forces preventing compression will deflect any incoming particle. This method, of course, violates balance of momentum, but this is plausible, since we can imagine that any momentum yielded by the fluid is grounded. More importantly though, dummy particles violate conservation of energy. To see this, consider a wave stopped by a wall. The fluid acts on the solid with equal and opposite force. Therefore, a wall particle $a$ is subject to acceleration, which should change its velocity after time step $\delta t$ to
	$$\bm{v}_a = \delta t \frac{\bm{F}_a}{m_a}.$$
	To preserve wall integrity, this velocity must be annulled, reducing the total energy by a small kinetic contribution. From this, we see that dummy particles are essentially dissipative. 
	
	One could argue that this is still physically reasonable because every fluid-solid collision creates a sound wave, leading to loss of energy. However, it is dubious whether such effect is correctly modeled by dummy particles. Moreover, it complicates efforts to make a strictly energy-conservative scheme.
	
	An alternative approach uses a fluid-solid force inspired by the Lennard-Jones potential 
	\begin{equation*}
		\bm{F}_{a \to b} =  -\frac{E_\text{wall}}{r} s_{ab}^2(1 - s_{ab}^2) \, \bm{e}_{ab}
	\end{equation*}
	by which a wall particle $a$ acts onto a fluid particle $b$. Here, 
	$$s_{ab} = \max \left\{ \frac{r_\text{wall}}{r_{ab}}, \; 1\right\}$$
	and $E_\text{wall}, r_\text{wall}$ are numerical parameters \cite{monaghan1994}. Except for this interaction, wall particles are neither considered in the balance of mass \eqref{bomo} nor in the calculation of density \eqref{boma-closed}. Unlike the dummy particle approach, this solid wall modeling is conservative, with a modified energy
	\begin{equation}
		\mathcal{H} = \sum_{a \in \text{fluid}} m_a\left( \frac{v_a^2}{2} + 
		\epsilon_a +  g z\right) + \sum_{a \in \text{wall}} \; \sum_{b \in \text{fluid}}E_\text{wall} \left( \frac{s_{ab}^2}{2} - \frac{s_{ab}^4}{4} - \frac{1}{4}\right).
		\label{boe+LJ}
	\end{equation}
	This makes it potentially advantageous for problems involving inviscid flows.
	
	In our case, when we seek a globally reversible-in-time simulation, the conservativeness of the fluid-wall interaction becomes important. However, even if we use the correct close formula for density \eqref{boma}, the initial state correction, and the conservative fluid-wall interaction, the symplectic Verlet scheme is not reversible globally in time because of the errors caused by the floating-point arithmetic. The following section contains the final ingredient necessary for that reversibility, the symplectic Verlet integrator with fixed-point arithmetic.
	
	\section{Symplectic integrator and fixed-point arithmetic}\label{sec.rev}
	The symplectic character of evolution equations \eqref{bomo} and \eqref{move} is preserved in symplectic numerical schemes, for instance the classical Verlet scheme \cite{leimkuehler}. Consequence are for instance the long-time stability of the trajectories and conservation of of integrals of motion (energy, momentum, and angular momentum, among others). However, even application of the Verlet scheme in SPH does not lead to globally in-time reversible simulations due to the floating-point errors, but a remedy is in the fixed-point arithmetic, as we show below.
	
\begin{figure}[ht!]
		\includegraphics[width=\columnwidth]{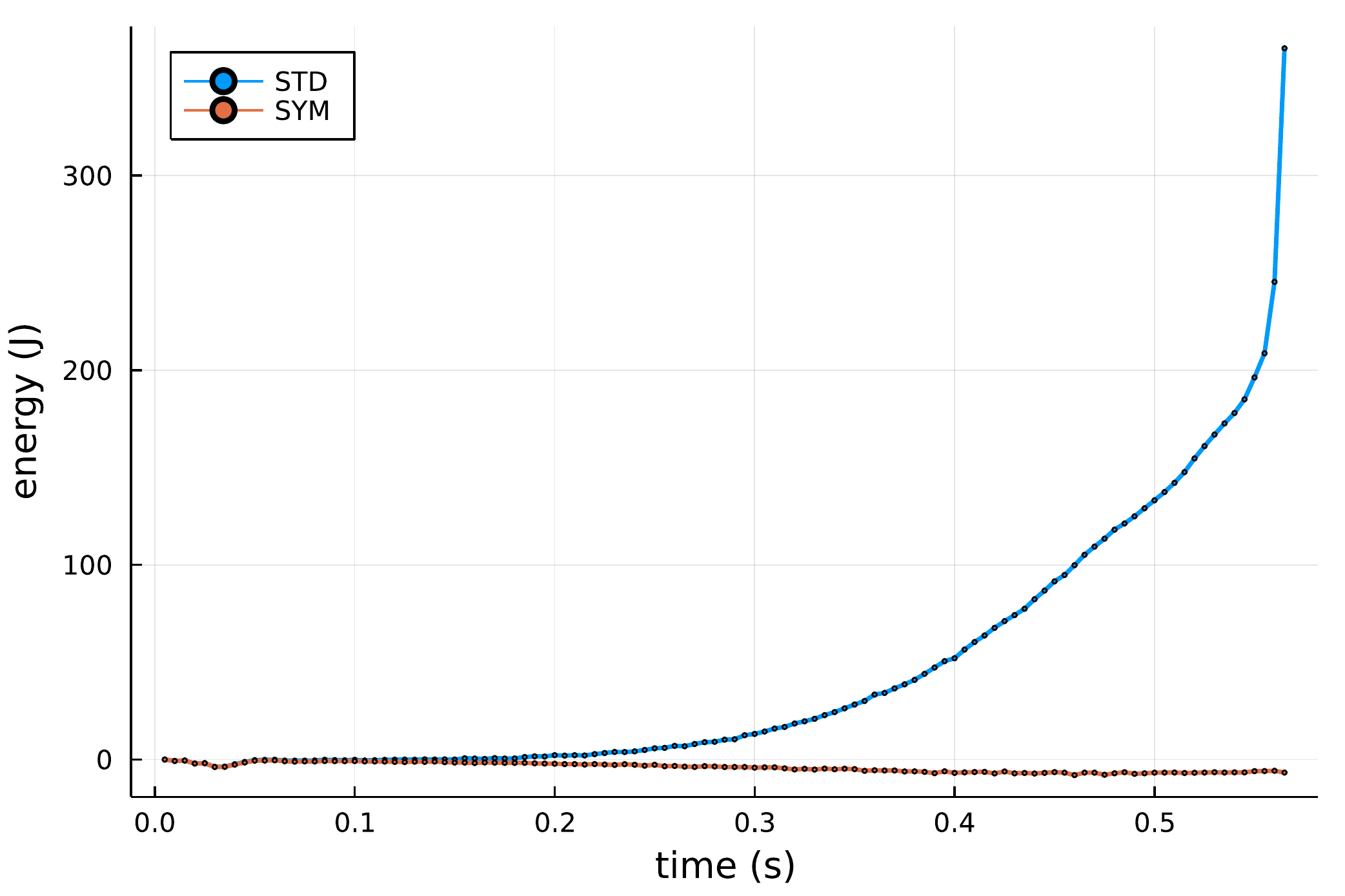}
		\caption{Energy growth in standard WCSPH (STD) and symplectic WCSPH (SYM) in dambreak simulation with zero viscosity and parameters from Table \ref{tab:DB_par}. Without a stabilization of some sort, STD diverges, whereas SYM is stable.}
		\label{fig:E_growth}
	\end{figure}

	The usual discretization of the SPH equations using the Verlet scheme is
	\begin{enumerate}
		\item $\bm{u}(t_{m+\frac{1}{2}}) = \bm{u}(t_{m}) + \tfrac{1}{2}\delta t \, \bm{a}(\bm{r}(t_m))$
		\item $\bm{r}(t_{m+1}) = \bm{r}(t_m) + \delta t \, \bm{u}(t_{m+\frac{1}{2}})$
		\item $\bm{u}(t_{m+1}) = \bm{u}(t_{m+\frac{1}{2}}) + \tfrac{1}{2}\delta t \, \bm{a}(\bm{r}(t_{m+1})),$
	\end{enumerate}
	where $\bm{a}$ is the acceleration composed of internal, gravitational, and wall forces. The scheme is of second order and symplectic, and therefore it conserves energy $\mathcal{H}_{\delta t}$ with error
	\begin{equation*}
		\mathcal{H}_{\delta t}(\bm{r}, \bm{u}) - \mathcal{H}(\bm{r}, \bm{u}) = O(\delta t^2)
	\end{equation*}
	that does not depend explicitly on time, provided that $\delta t$ is sufficiently small \cite{candy1991}\footnote{Here, we rely on the fact that $p = p(\rho) = p(\rho(\bm{r}))$ using the closed expression \eqref{boma-closed}.}. If we took $\rho$ as a separate state variable with its own evolution equation \eqref{boma} (usual WCSPH),  symplecticity of the system of equations would be violated. The usual WCSPH equations are still Hamiltonian (generated by a Poisson bracket and a Hamiltonian), as checked by program \cite{kroeger2010}, but they are not symplectic and they thus the (symplectic) Verlet scheme does not preserve their geometric structure, which causes energy to grow exponentially, unless some stabilization is used - see Figure  \ref{fig:E_growth}.\footnote{In \cite{monaghan2005}, it is suggested to correct this error by a mid-step extrapolation of density. However, it should be noted this approach adds additional inaccuracy to the scheme.}
	
	Besides being symplectic, Verlet scheme is also time reversible, which means that after changing the sign of velocities, the scheme returns back in time. The scheme thus preserves the time-reversal symmetry of the original system \eqref{boma}, \eqref{bomo}, \eqref{move}. We can verify this by substituting 
	\begin{align*}
		&\bm{u}(t_{m}) \mapsto -\bm{u}(t_{m+1})\\
		&\bm{u}(t_{m+1}) \mapsto -\bm{u}(t_{m})\\
		&\bm{r}(t_{m}) \mapsto \bm{r}(t_{m+1})\\
		&\bm{r}(t_{m+1}) \mapsto \bm{r}(t_{m})
	\end{align*}
	which produces the same set of equations
	\begin{enumerate}
		\item $-\bm{u}(t_{m+\frac{1}{2}}) = -\bm{u}(t_{m+1}) + \tfrac{1}{2}\delta t \, \bm{a}(\bm{r}(t_{m+1}))$
		\item $\bm{r}(t_{m}) = \bm{r}(t_{m+1}) - \delta t \, \bm{u}(t_{m+\frac{1}{2}})$
		\item $-\bm{u}(t_{m}) = -\bm{u}(t_{m+\frac{1}{2}}) + \tfrac{1}{2}\delta t \, \bm{a}(\bm{r}(t_{m})),$
	\end{enumerate}
	but in the reversed order. Unfortunately, the time reversibility fails in the floating-point arithmetic (FloPA), where addition is \emph{not} associative. In particular, addition and subtraction of a value $\delta x$ to a float $x$ does not recover $x$. For instance
	\begin{equation}
		(1 + 0.5\epsilon) - 0.5\epsilon \overset{\text{float}}{=} 0.9999999999999999 \neq 1.
		\label{nonassoc}
	\end{equation}
	Although this error is usually very small, it tends to accumulate in the presence of non-linearities and destroys the reversibility in longer simulations. 
	
\begin{figure}[ht!]
		\begin{subfigure}{\textwidth}
			\includegraphics[trim = {0cm 0cm 0cm -4cm}, width=\linewidth]{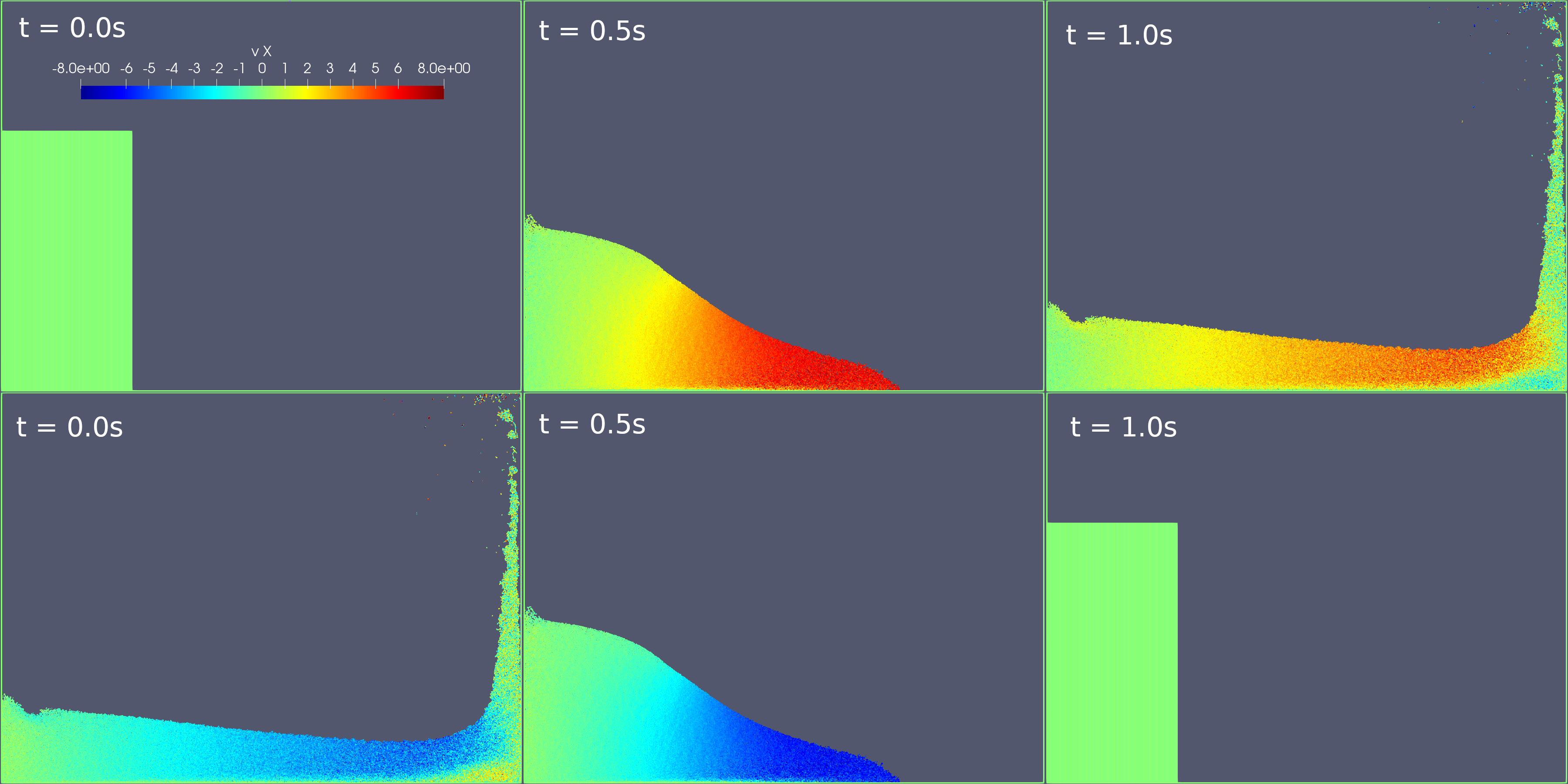}
			\caption{FixPA}
		\end{subfigure}	
		\begin{subfigure}{\textwidth}
			\includegraphics[width=\linewidth]{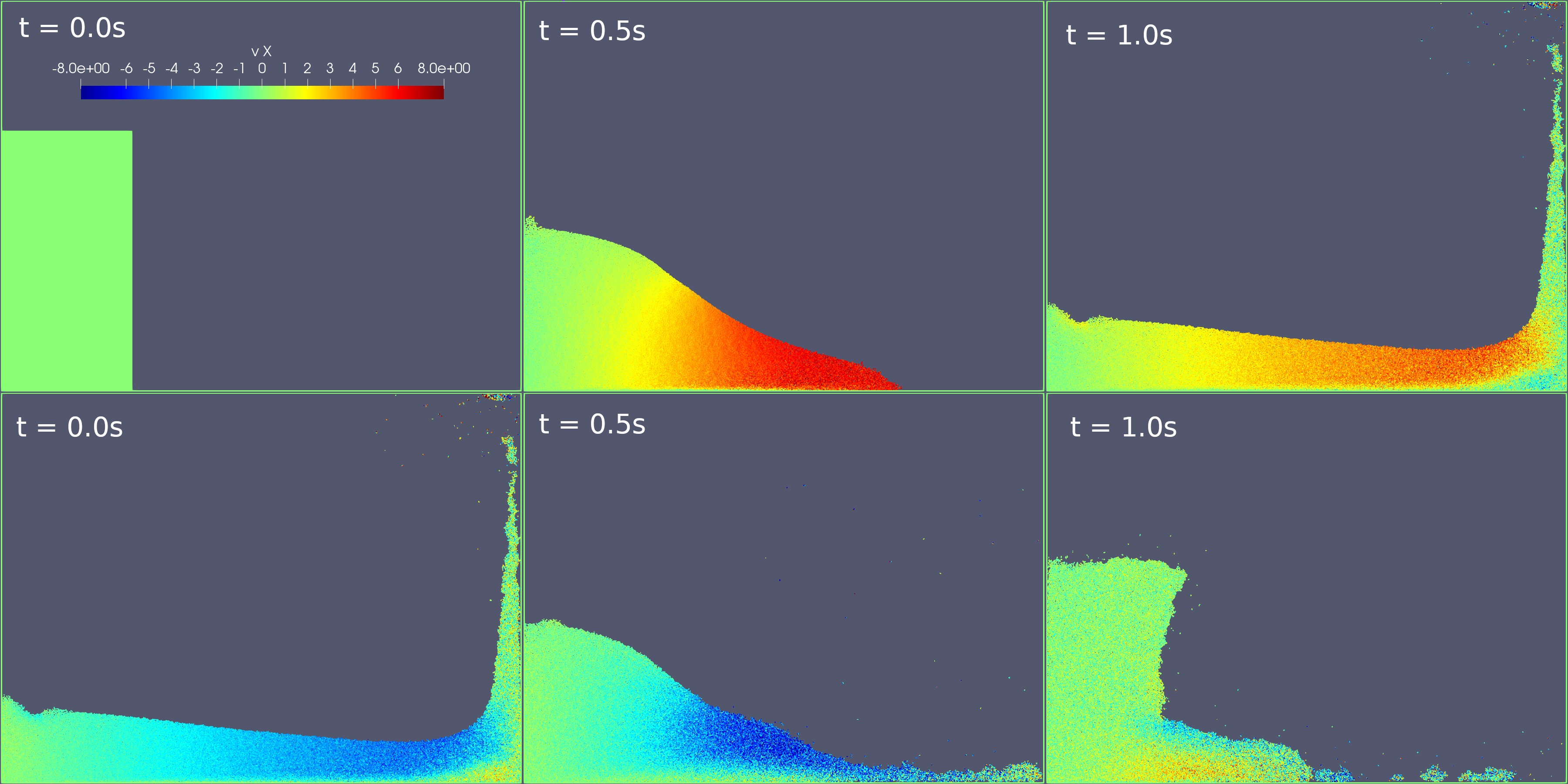}
			\caption{FloPA}
		\end{subfigure}
		\caption{Results for dam break test with parameters from Table~\ref{tab:DB_par} using symplectic WCSPH. At $t = 1\text{s}$, we reset the time counter and revert all velocities. In FixPA, the simulation is completely reversible and returns to the initial state. In FloPA, this fails due to round-off errors illustrated in \eqref{nonassoc}.}
		\label{FixPa_vs_FloPa}
	\end{figure}

	There is, however, an easy solution suggested by \cite{JANUS} in the context of N-body simulations of Solar System. It converts vectors $\bm{r}(t_m)$, $\bm{u}(t_m)$, and $\bm{a}(t_m)$ to the fixed-point arithmetic (FixPA) just before they are used in the time-step evaluation. Since addition is associative in FixPA, this method allows for bit-precise time reversibility and long-time conservation of energy -- see Figure \ref{FixPa_vs_FloPa}. Note that the intermediate computations of the sums in \eqref{boma-closed} and \eqref{bomo} can be still performed in FloPA and only then converted to FixPA. However, due to the non-associativity of FloPA, summation of more than two elements in FloPA is order-dependent, and to ensure reversibility, the order of the summands must be chosen in a way that does not depend on the current state of the simulation (for instance, it is possible to perform summations always in the order of the particle indices).
	
	\begin{figure}[ht!]
		\centering
		\includegraphics[width=\textwidth]{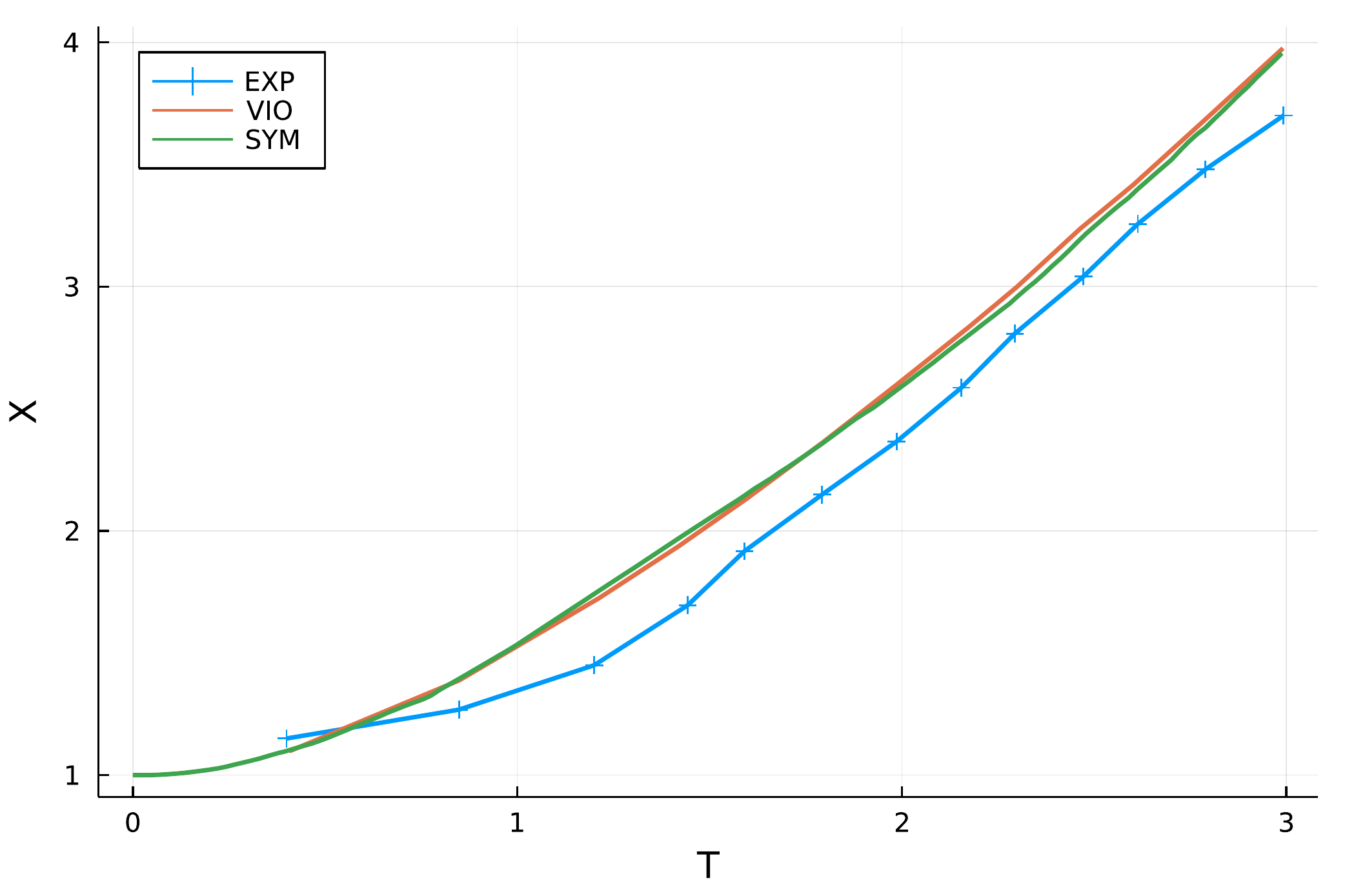}
		\caption{Comparison of symplectic WCSPH (SYM) in dam-break simulation with experimental data \cite{koshizuka96} (EXP) and numerical result from Violeau's book \cite{violeau} (VIO). Here, $T = t\sqrt{g/l_\text{wcv}}$ is a dimensionless time and $X = l_\text{le}/l_\text{wcv}$, where $l_\text{le}$ is the $x$-coordinate of the wave's leading edge.}
	\end{figure}
	
	\begin{table}[ht!]
		\centering
		\begin{tabular}{|l|c|c|}
			\hline
			density & $\rho$ & 1000 kg/$\text{m}^3$\\
			\hline
			spatial step & $\dd{r}$ & 0.005
			$\text{m}$\\
			\hline
			num. sound speed & $c$ & 120 m/s\\
			\hline
			gravitational constant & $g$ & 9.8 m/$\text{s}^2$\\
			\hline
			water column width & $l_\text{wcw}$ & 1m \\
			\hline
			water column height & $l_\text{wch}$ & 2m \\
			\hline
			box width & $l_\text{bw}$ & 4m \\
			\hline
			box height & $l_\text{bh}$ & 3m \\
			\hline
			kernel support radius & $h$ & $3\dd{r}$\\
			\hline
			Lennard-Jones radius & $r_\text{wall}$ & $0.95\dd{r}$\\
			\hline
			Lennard-Jones energy & $E_\text{wall}$ & $10 \, m g l_\text{wch}$\\
			\hline
			time step & $\dd{t}$ & $0.2 \, h/c$\\
			\hline
		\end{tabular}
		\caption{Parameters used for dam-break test.}
		\label{tab:DB_par}
	\end{table}
	
	In summary, symplectic WCSPH with the correct closed formula for density, the initial state correction, the conservative wall-fluid interaction, and the symplectic numerical Verlet scheme in the fixed-point arithmetic finally lead to globally-in-time reversible simulations. 
	We have demonstrated the reversibility for instance on the dam break benchmark, and Appendix \ref{sec.Gresho} contains the Gresho vortex benchmark.
	Although the simulation is symmetric, we can still observe the growth of entropy if we choose not see all the details of the simulation, which is the matter of the following section.
	
	\section{Emergence of the second law of thermodynamics}

	The second law of thermodynamics tells that the entropy of each isolated system grows until the system reaches the thermodynamic equilibrium. Where does this irreversible behavior come from when the underlying evolution equations for particle dynamics (here SPH) are reversible? In this section we illustrate that the emergence of the second law from completely reversible dynamics is caused by ignoring details of the dynamics, similarly as in \cite{JSP2020}.
	
\begin{figure}[ht!]
		\includegraphics[width=0.5\columnwidth]{"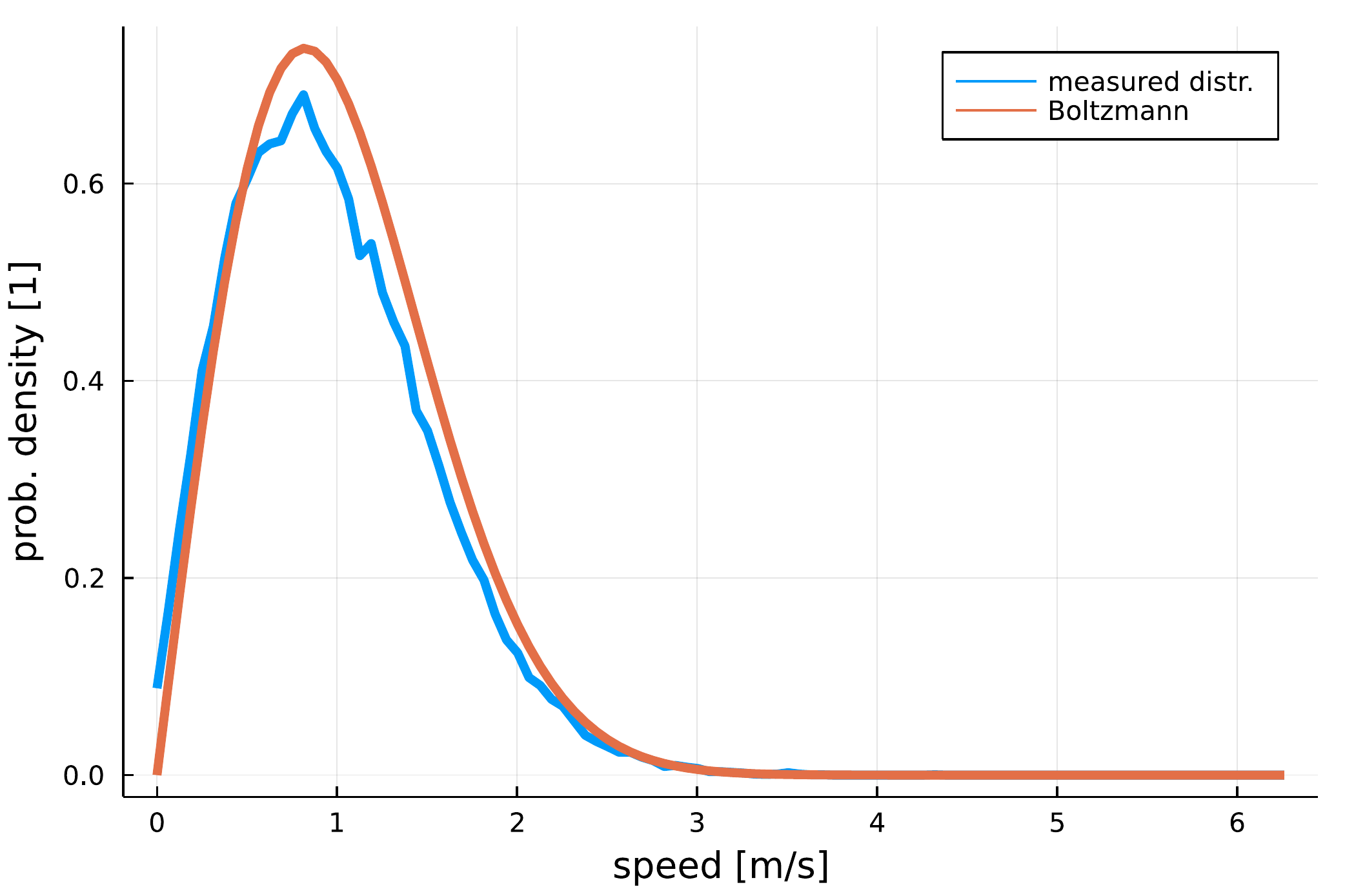"}
		\includegraphics[width=0.5\columnwidth]{"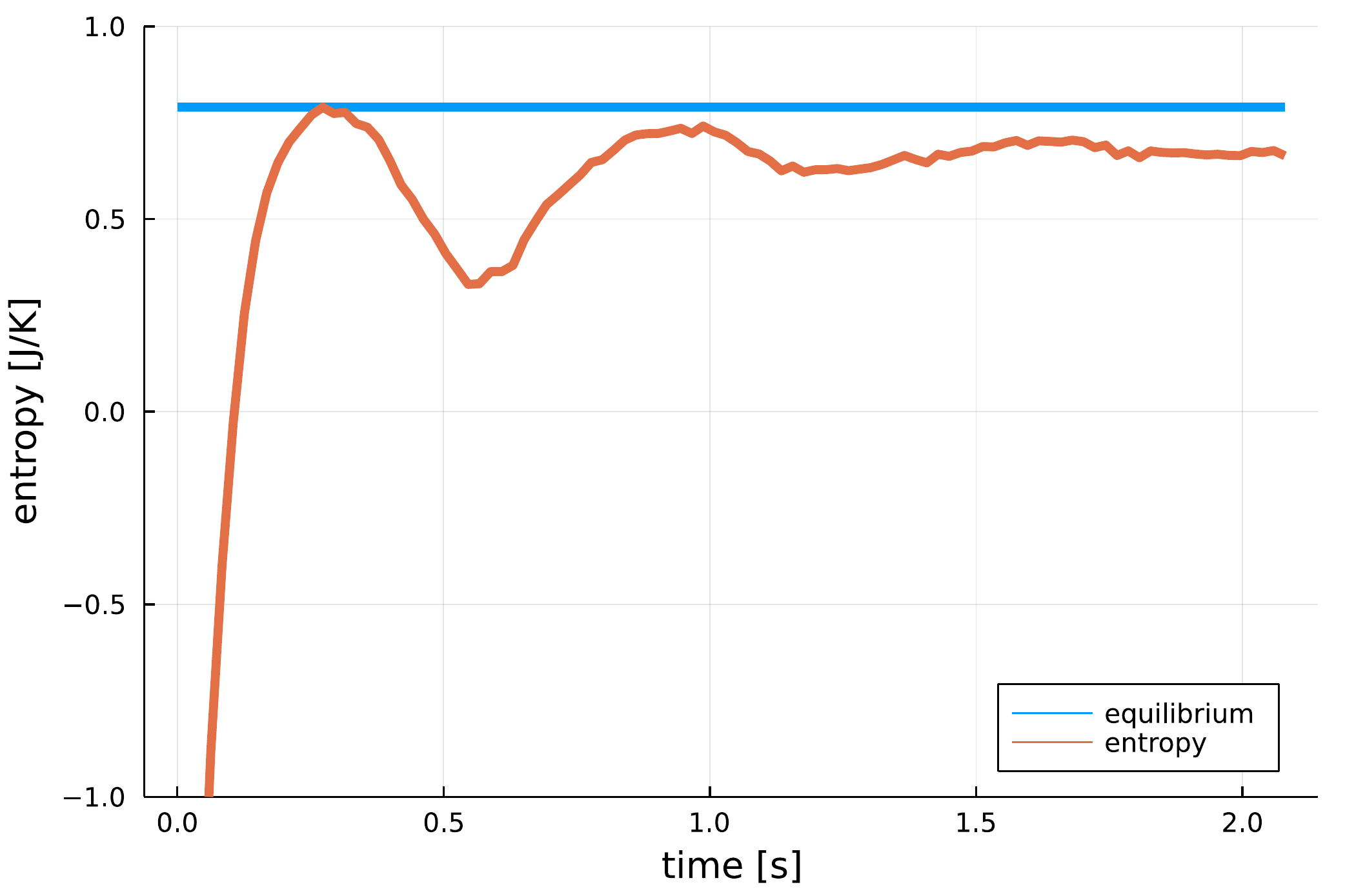"}
		\caption{\label{fig.entropy} 
			\textit{Left:} Histogram of the velocities of the particles at time $t=2.0$ (blue) in comparison with the equilibrium Maxwell-Boltzmann distribution (orange) \eqref{eq.MB.eq} with fitted temperature $T$.
			\textit{Right:} Evolution of Boltzmann entropy (see Appendix \ref{sec.Boltzmann}) in the reversible breaking-dam simulation (blue). Despite the reversibility of the dynamics, the entropy increases and approaches equilibrium entropy corresponding to the temperature of the measured Maxwell-Boltzmann distribution (blue line) from \eqref{eq.S.eq.E}.}
	\end{figure}

	The emergence of the second law is actually expected due to the result of Lanford \cite{lanford} and following works \cite{pulvirenti}, where a system of classical particles with a short-range potential is shown to obey the Boltzmann equation. The collision term in the Boltzmann equation then causes the growth of the Boltzmann entropy. Figure \ref{fig.entropy} shows the growth of Boltzmann entropy in a reversible SPH breaking-dam simulation.
	But how to interpret the Lanford's mathematical result and the observed growth of Boltzmann entropy of the SPH particles from the physical point of view? The SPH particles obey reversible Hamiltonian dynamics. Another possibility to describe their motion is the Liouville equation, which expresses reversible evolution of the N-particle distribution function, $f_N(t, \rr_1, \pp_1, \dots, \rr_N, \pp_N)$. Indeed, Liouville equation follows from the particle mechanics, and, vice versa, if we set the distribution function to be the product of Dirac $\delta$-distributions in the positions and momenta of all particles, we recover Hamilton canonical equations. Liouville equation also conserves the Liouville entropy 
	\begin{equation}
		S^{\text{Liouville}}(f_N) = -\frac{k_B}{N!}\int d\rr_1 \int d\pp_1\dots \int d\rr_N \int d\pp_N f_N \ln\left(h^{3N}f_N\right),
	\end{equation}
	which thus remains constant. This is consistent with the reversibility of the underlying Hamiltonian particle mechanics.

	Now we decide not to see all the positions and momenta of the individual particles, observing for instance only the probability distribution of momenta at given space and time. Such one-particle distribution function $f(t,\rr,\pp)$ indeed does not contain the knowledge of positions and momenta of all particles and thus the entropy of the system expressed in terms of $f$ must be higher than the Liouville entropy (functional of $f_N$). The former is the Boltzmann entropy and it is obtained by maximization of the Liouville entropy subject to the constraint that the one-particle distribution function is known \cite{pkg}.
	Plugging the resulting $f_N$ back into the Liouville entropy, we obtain the Boltzmann entropy 
	\begin{equation}
		S^{\text{Boltzmann}}(f) = -k_B\int d\rr d\pp f (\ln (h^2 f)-1)
	\end{equation}
	see \cite{pkg}. Because the Boltzmann entropy is the Liouville entropy evaluated at the point $f_N$ where it is maximal (keeping the knowledge of $f$), the Boltzmann entropy is higher than the Liouville entropy. Although the Liouville entropy remains constant in the dynamics, the Boltzmann entropy can grow in the dynamics following the Liouville equation.
	
	The irreversibility of the dynamics of the one-particle distribution function can be illustrated on the reversible breaking-dam simulation with fixed-point arithmetic. Section \ref{sec.rev} contains details of the simulation, but let us also show the histogram of the particle velocity distribution in the middle of the simulation, just before the velocities are inverted. Figure \ref{fig.entropy} shows the measured histogram of the velocities of the SPH particles in comparison with the equilibrium Maxwell-Boltzmann distribution function \eqref{eq.MB.eq} obtained by fitting the effective temperature\footnote{Note that the effective temperature is not the physical temperature of the system because it is calculated from the overall motion of the macroscopic SPH particles. To see this, consider the extreme case of only one SPH particle with weight one kilogram, where the effective temperature obtained from the motion of that macroscopic particle does not coincide with the physical temperature of that particle.}. Although all velocities were initially zero, their distribution approaches the equilibrium distribution. This approach towards equilibrium is reflected in the growth of Boltzmann entropy during the long-time reversible simulation.
	
	In other words, when the particles are described by the Boltzmann equation, our knowledge is incomplete, since the positions and momenta of all the particles remain unknown, in contrast to the description by means of Hamilton canonical equations or Liouville equation. And this lack of knowledge makes the evolution of the one-particle distribution function $f$ irreversible. The irreversibility is then explicitly expressed by the collision integral in the Boltzmann equation. The second law of thermodynamics thus emerges from completely reversible dynamics when our description is incomplete (not seeing all positions and momenta of the particles).

	In summary, if we see all the positions and momenta in the SPH simulation, we can not see the second law of thermodynamics. Indeed, the simulation is reversible and the Liouville entropy remains constant. However, when we only focus on the one-particle distribution function, we can see the growth of Boltzmann entropy and thus irreversible behavior. The second law of thermodynamics emerges when we our knowledge about the precise state of the system is incomplete \cite{JSP2020}.
	
	Interestingly, since our discrete system is deterministic, reversible and can exist in only finitely many states, it is recurrent \cite{wallace2015recurrence}. In other words, if the system evolves long enough, it comes back to the initial state exactly (provided that the particles are not allowed to escape to infinity). However, since the phase space of the system contains thousands of SPH particles, it has enormous amount of states, which makes it highly unlikely to observe the recurrence theorem in practice.
	
	\section{Conclusion}
	
	In this paper we have turned the usual weakly compressible smoothed particle hydrodynamics (WCSPH) to a symplectic form by finding the correct closed formula for the mass density. Because the closed formula for density leads to different treatment of particles that are initially near the boundary, the method of initial state correction (ICS) was introduced, due to which all particles are then treated in the same way. This leads to stable SPH simulations in the presence of free surfaces without any further stabilization. 
	
	In order to get simulations that preserve the energy, we use a conservative fluid-wall interaction. A symplectic (Verlet) scheme, which is suitable for the symplectic WCSPH, implemented in the fixed-point arithmetic then leads to globally-in-time reversible SPH simulations. This is demonstrated on the dam break benchmark, where inversion of velocities at a later stage of the simulation eventually leads the system back to the initial state. The simulations are available in a new Julia package \textit{SmoothedParticles.jl} \cite{Kincl_SmoothedParticles_jl}.
	
	Despite the reversibility of the simulations, we observe thermodynamic behavior when we do not use all the details of the simulation. The Boltzmann entropy, which depends only on the one-particle distribution function and not on positions of individual particles, grows in the dam-break simulation and approaches the equilibrium value. In other words, we observe the emergence of the second law of thermodynamics from purely reversible dynamics, caused by reduction of our knowledge about the system.
	
	In future, we would like to extend the WCSPH framework to non-isothermal fluids and solids while keeping the Hamiltonianity of the equations and numerical schemes.

	\appendix
	\section{Gresho vortex benchmark}\label{sec.Gresho}
	The Gresho vortex benchmark \cite{grimmstrele-2014} prescribes a tangential velocity component of Eulerian fluid in polar coordinates as
	\begin{equation*}
		u^0_\theta(r) = \begin{cases}
			5r & \mathrm{for} \quad r < \frac{1}{5},\\
			2 - 5r & \mathrm{for} \quad \frac{1}{5} < r < \frac{2}{5},\\
			0 & \mathrm{for} \quad r > \frac{2}{5}.
		\end{cases}
	\end{equation*}
	The vortex is confined in a box $(-\frac{1}{2},\frac{1}{2}) \times (-\frac{1}{2},\frac{1}{2})$ with no-slip wall (though different variants can be found in literature). Theoretically, in the absence of viscosity, the vortex should be stationary, but obtaining this in a numerical scheme is challenging. In our numerical experiment, we simulate the flow in time-interval $\left[0,1\right]$ and measure the error using a discretized $(L^\infty,L^2)$ Sobolev norm
	\begin{equation}
		e = \max_{t \in \left[0,1\right]} \sqrt{\frac{75}{4\pi}\sum_{a } V_a \left| \bm{u}_a - \bm{u}^0(\bm{r_a})\right|^2}.
		\label{gresho-error}
	\end{equation}
	The meaning of the factor $\frac{75}{4\pi}$ is that the zero velocity field corresponds to error $1$. In order to prevent formation of void space in the center of vortex, we had to include an anti-clump force term in the form
	\begin{figure}[h]
		\centering
		\includegraphics[width=\linewidth]{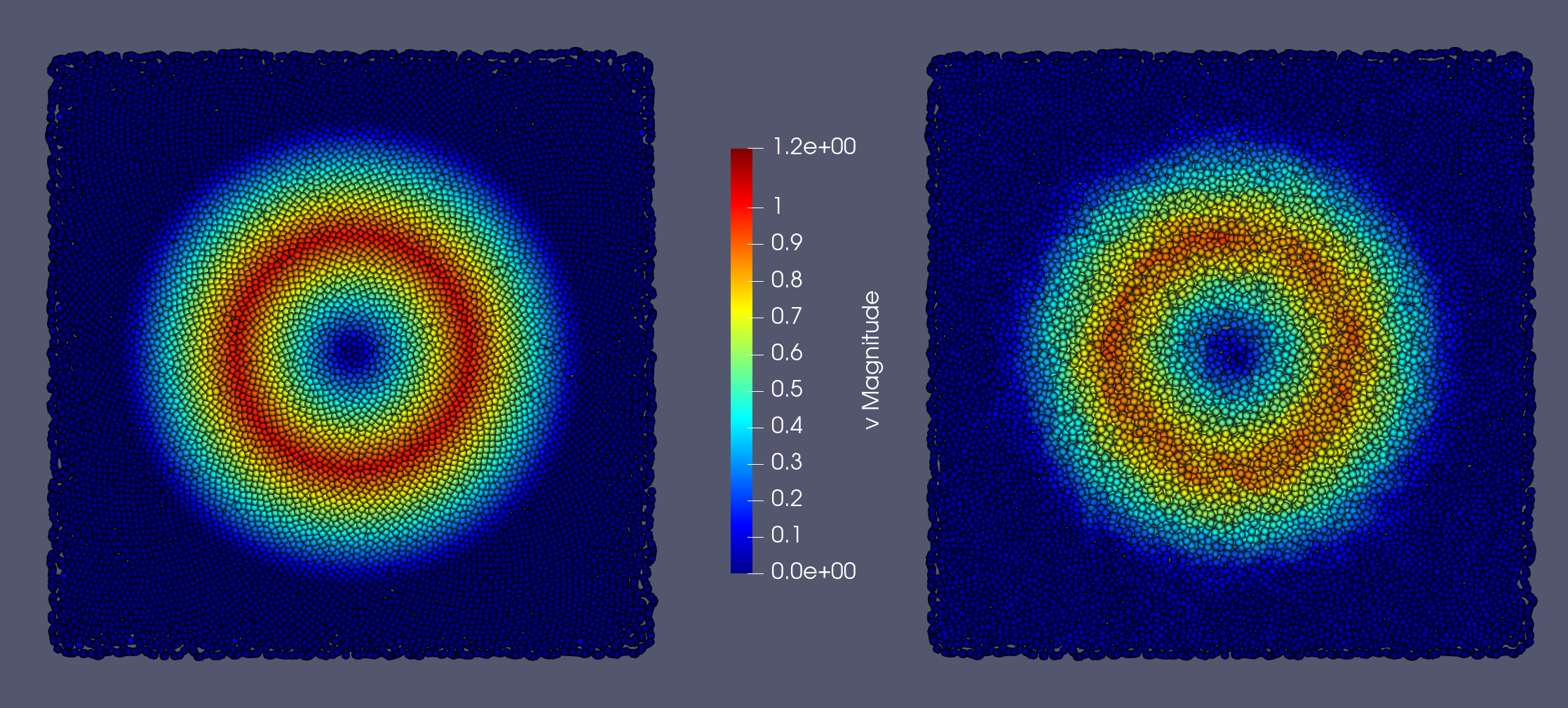}
		\caption{Simulation result with Vogel spiral + ISC arrangement at $t = 0$ (left) and $t = 1$ (right).}
		\label{fig.gresho}
	\end{figure}
	\begin{table}[]
		\centering
		\begin{tabular}{|l|l|l|l|l}
			\cline{1-4}
			& no filter & passive filter & active filter &  \\ \cline{1-4}
			square      & 28.16\%   & 12.20\%        & 53.07\%       &  \\ \cline{1-4}
			hexagonal   & 27.84\%   & 11.77\%        & 55.79\%       &  \\ \cline{1-4}
			Vogel       & 29.85\%   & 13.31\%        & 78.06\%       &  \\ \cline{1-4}
			Vogel + ISC & 28.19\%   & 11.26\%        & 55.25\%       &  \\ \cline{1-4}
		\end{tabular}
		\caption{Error of Gresho vortex benchmark for different grid types (square, hexagonal, Vogel spiral, Vogel spiral + ISC) and with different types of noise filtering. The error is computed according to formula \eqref{gresho-error}. It appears from data that passive (a posteriori) filter is better than active filtering (which dissipates energy). Compared to WENO method \cite{grimmstrele-2014}, numerical dissipation in symplectic WCSPH is quite strong.}
		\label{tab.gresho-error}
	\end{table}
	\begin{table}
		\centering
		\begin{tabular}{|l|c|c|}
			\hline
			density & $\rho$ & 1\\
			\hline
			spatial step & $\dd{r}$ & 1e-2\\
			\hline
			num. sound speed & $c$ & 20\\
			\hline
			kernel support radius & $h$ & $3\dd{r}$\\
			\hline
			box size & $l$ & 1\\
			\hline
			anti-clump pressure & $p_0$ & 10\\
			\hline
			time step & $\dd{t}$ & $0.1h/c$\\
			\hline
			filter frequency & $M$ & 30 \\
			\hline
		\end{tabular}
		\caption{Paramers used for Gresho vortex simulation. This benchmark is dimensionless and with no gravity. Walls were implemented by two layers of dummy particles.}
		\label{tab.gresho-par}
	\end{table}
	\begin{equation*}
		\bm{f} = \sum_b m_b \left(\frac{p_0}{\rho_a^2} + \frac{p_0}{\rho_b^2}\right) \tilde{w}'(r_{ab}) \bm{e}_{ab},
	\end{equation*}
	where $\tilde{w}$ is a smoothing kernel with support radius $\frac{\dd{r}}{2}$. Additionally, we consider noise reduction by the Shepard filter
	\begin{equation*}
		\bm{\tilde{u}}_a = \frac{1}{\gamma_a}\sum_b V_b \bm{u}_b w_{ab},
	\end{equation*} 
	where
	\begin{equation*}
		\gamma_a = \sum_b V_b w_{ab}.
	\end{equation*}
	This filter can be applied either in post-processing (passive filter), or by setting $\bm{u}_a := \bm{\tilde{u}}_a$ every $M$ time steps (active filter). Results are summarized in Table \ref{tab.gresho-error} and Figure \ref{fig.gresho}. Simulation parameters are listed in Table \ref{tab.gresho-par}.

	\section{The Maxwell-Boltzmann entropy}\label{sec.Boltzmann}
	In this section we derive the formula for Boltzmann entropy in terms of the Maxwell-Boltzmann distribution and, subsequently, to find its equilibrium value. 
	
	\subsection{The reduced Boltzmann entropy in terms of the Maxwell-Boltzmann distribution}
	Let us start with the Boltzmann entropy in two dimensions expressed in terms of the one-particle distribution function $f(t,\rr,\pp)$,
	\begin{equation}
		S^{(\text{Boltzmann})} = -k_B \int d\rr\int d\pp f(\ln(h^2 f) -1),
	\end{equation}
	where the position is constrained to a box with volume $V$ and where $h$ is the Planck constant.
	This formula can be obtained for instance by the principle of maximum entropy from the Liouville entropy, from where it also follows that the one-particle distribution function is normalized to the number of particles, see \cite{pkg}. 
	
	Assuming that the distribution function depends only on the magnitude of the momentum (isotropic dependence), the Boltzmann entropy can be rewritten as
	\begin{equation}
		S^{(\text{Boltzmann})} = k_B N\ln \frac{e}{h^2} -2\pi k_B V \int_0^\infty d p\, p f(p)\ln f(p),
	\end{equation}
	where $p = |\pp|$ is the norm of the momentum. This dependence on the one-particle distribution function $f(p)$ has to be converted to a dependence on the Maxwell-Boltzmann distribution function $f_{MB}(v)$, which tells the probability that the norm of velocity of a particle $v=p/m$ is in the interval $(v,v+dv)$. This is done using the normalization, 
	\begin{equation}
		1= \frac{1}{N} \int d\rr\int d\pp f = \int_0^\infty dv\, \frac{2\pi m^2 V}{N} v f(v),
	\end{equation}
	which leads to the expression of the two-dimensional Maxwell-Boltzmann distribution function in terms of the one-particle distribution,
	\begin{equation}
		f_{MB}(v) = \frac{2\pi m^2 V}{N} v f(v),
	\end{equation}
	which is normalized to unity. Note that we assume that the one-particle distribution is homogeneous in space (the total volume being $V$) and isotropic in momentum. These assumptions are valid in the thermodynamic equilibrium, but they approximately hold also before the equilibrium is reached, and for the purpose of showing that the Boltzmann entropy grow in our simulations, the approximation is satisfactory. Finally, the Boltzmann entropy in terms of the Maxwell-Boltzmann distribution function becomes
	\begin{equation}
		S^{(\text{Boltzmann})} - k_B N\ln \frac{e}{h^2}= -k_B N \int_0^\infty d v\, f_{MB}(v)\ln \left(\frac{f_{MB}(v)N}{2\pi m^2 v}\right).
	\end{equation}
	Since the number of particles is conserved in our simulation, we only evaluate the part of Boltzmann entropy that varies when $f_{MB}$ changes, 
	\begin{align}
		S^{(\text{Boltzmann})}_{\text{reduced}} &= \frac{S^{(\text{Boltzmann})} - k_B N\ln \frac{e}{h^2}}{k_B N} + \ln\left(\frac{N}{2\pi m^2}\right)\nonumber\\
		&=-\int_0^\infty f_{MB}(v) \ln \left(\frac{f_{MB}(v)}{v}\right),
	\end{align}
	called reduced Boltzmann entropy. Note that the integral converges near the origin because the Maxwell-Boltzmann distribution typically has linear growth there. 
	
	In our simulations, this reduced Boltzmann entropy is numerically integrated using the approximate Maxwell-Boltzmann distribution function obtained from a histogram of the particle velocities at given time instant.
	
	\subsection{Equilibrium Boltzmann entropy}
	What is the final value of the reduced Boltzmann entropy when the system of particles reaches the thermodynamic equilibrium? This question can be answered in two steps. First, the equilibrium distribution function is calculated by the principle of maximum entropy (MaxEnt). Second, the equilibrium distribution function is plugged into the formula for the reduced Boltzmann entropy. Note, however, that the MaxEnt step depends on the energy (Hamiltonian) of the system and it depends on the complexity of the Hamiltonian whether the calculation can proceed purely analytically (without numerical solutions). Therefore, we use the simplest Hamiltonian approximating the true Hamiltonian of our SPH particles, consisting only of the kinetic energy of the particles. Maximization of the Boltzmann entropy subject to the constraints given by the total energy and total number of particles leads to the equilibrium one-particle distribution function
	\begin{equation}
		f_{\text{equilibrium}} = \frac{1}{h^2} e^{-\frac{N^*}{k_B}} e^{-\frac{E^*}{k_B}\frac{\pp^2}{2m}},
	\end{equation}
	where $N^*$ and $E^*$ are the Lagrange multipliers corresponding to the two constraints (number of particles and total energy). From the normalization it follows that $N/V = \exp(-N^*/k_B)/h^2$, while the other Lagrange multiplier can be interpreted as the inverse temperature, $E^* = 1/T$. 
	
	The equilibrium Maxwell-Boltzmann distribution function then becomes
	\begin{equation}\label{eq.MB.eq}
		f_{MB,\text{equilibrium}}(T,v) = \frac{m}{k_B T}v e^{-\frac{1}{2}m v^2}.
	\end{equation}
	As time proceeds in our simulations, the histogram of particle velocities approaches the equilibrium Maxwell-Boltzmann distribution.
	
	When this equilibrium distribution function is plugged back into the Boltzmann entropy, the reduced Boltzmann entropy becomes
	\begin{equation}\label{eq.S.eq.T}
		S^{(\text{Boltzmann})}_{(\text{reduced})}(T) = 1 + \ln\left(\frac{k_B T}{m}\right).
	\end{equation}
	This equilibrium entropy, which depends on $T$, can be calculated once the temperature is obtained by fitting the histogram of particle velocities to the equilibrium Maxwell-Boltzmann distribution function \eqref{eq.MB.eq}. This makes sense, however, only at later stages in the simulations, when the histogram approaches the equilibrium distribution.
	
	Instead of using the equilibrium temperature, which makes sense only in later stages of the simulations, we can express the equilibrium entropy in terms of energy, which can be measured anytime. The total kinetic energy of the two-dimensional system of particles is equal to $E = N k_B T$, which makes it possible to write the equilibrium entropy in terms of the energy,
	\begin{equation}\label{eq.S.eq.E}
		S^{(\text{Boltzmann})}_{(\text{reduced})}(E) = 1 + \ln\left(\frac{E}{Nm}\right).
	\end{equation}
	This value of equilibrium entropy is close to the Boltzmann entropy obtained directly from the approximated Maxwell-Boltzmann distribution, while the entropy based on the equilibrium temperature is only approached at later stages of the simulations. 
	
	\section*{Acknowledgments}
	OK was supported by project No. START/SCI/053 of Charles University Research program.
	MP was supported by project No. UNCE/SCI/023 of Charles University Research program. OK and MP were also supported by the Czech Science Foundation (project no. 20-22092S)
	
\end{document}